\documentclass[12pt,a4paper]{amsart}

\usepackage[T1]{fontenc}
\usepackage[utf8]{inputenc}
\usepackage[british]{babel}
\usepackage{mathtools}
\usepackage{amsthm}
\usepackage{libertine}
\usepackage[libertine]{newtxmath}
\usepackage[mathscr]{euscript}
\usepackage{enumitem}
\usepackage{tikz-cd}
\usetikzlibrary{decorations.markings}
\usepackage{float}
\usepackage[
  backend=biber,
  style=alphabetic,
  maxnames=10,
  maxalphanames=10]{biblatex}
\usepackage{hyperref}
\usepackage[noabbrev]{cleveref}

\theoremstyle{plain}
\newtheorem{thm}{Theorem}
\newtheorem*{thm*}{Theorem}
\newtheorem{lm}[thm]{Lemma}
\newtheorem{defnlm}[thm]{Definition-Lemma}

\newtheorem{cor}[thm]{Corollary}

\theoremstyle{definition}
\newtheorem{defn}[thm]{Definition}

\newtheorem{nota}[thm]{Notation}
\newtheorem{exmp}[thm]{Example}

\theoremstyle{remark}
\newtheorem{rem}[thm]{Remark}
\Crefname{thm}{Theorem}{Theorems}
\Crefname{lm}{Lemma}{Lemmas}
\Crefname{prop}{Proposition}{Propositions}
\Crefname{cor}{Corollary}{Corollaries}
\Crefname{fact}{Fact}{Facts}
\Crefname{q}{Question}{Questions}
\Crefname{defn}{Definition}{Definitions}
\Crefname{construction}{Construction}{Constructions}
\Crefname{nota}{Notation}{Notations}
\Crefname{exmp}{Example}{Examples}
\Crefname{xca}{Exercise}{Exercises}
\Crefname{rem}{Remark}{Remarks}

\title[Extension properties of adapted differentials on klt orbifolds]{Extension of adapted differentials on klt orbifolds}
\author[Pedro N{\'u}{\~n}ez]{Pedro N{\'u}{\~n}ez}
\address{Pedro N\'{u}\~{n}ez \newline
\indent National Taiwan University \newline
\indent Department of Mathematics \newline
\indent No.~1, Sec.~4, Roosevelt Rd., \newline
\indent Taipei 10617, Taiwan}
\email{\normalfont\href{mailto:pnunez@ntu.edu.tw}{pnunez@ntu.edu.tw}}

\urladdr{\normalfont\href{https://homepage.ntu.edu.tw/~pnunez}{https://homepage.ntu.edu.tw/\textasciitilde pnunez}}
\thanks{The author was supported by the DFG-Graduiertenkolleg GK1821 and by the NSTC of Taiwan, with grant numbers 112-2811-M-002-108 and 113-2811-M-002-087.}

\date{\today}

\makeatletter
\hypersetup{
  pdfauthor={\authors},
  pdftitle={\@title},
  colorlinks,
  linkcolor=[rgb]{0.2,0.2,0.6},
  citecolor=[rgb]{0.2,0.2,0.6},
  urlcolor=[rgb]{0.2,0.2,0.6}}
\makeatother

\begin{document}

\maketitle

\begin{abstract}
  Given a geometric orbifold $(X,\Delta)$ in the sense of Campana,
  adapted reflexive differentials with respect to this orbifold
  are defined on suitably ramified covers of $X$.
  We show that if the orbifold $(X,\Delta)$ is klt,
  then any such reflexive differential form can be extended
  to a regular differential form on a resolution of singularities of the cover.
\end{abstract}

\tableofcontents

\section{Introduction}

\subsection{Motivation}

Geometric orbifolds, also known as \emph{$\mathcal{C}$-pairs},
were introduced by Campana in \cite{cam04} as part of his program
for the birational classification of algebraic varieties.
They are pairs $(X,\Delta)$ consisting of a normal variety $X$
and a Weil $\mathbb{Q}$-divisor of the form $\sum_{i} \frac{m_{i}-1}{m_{i}}D_{i}$
for $m_{i} \in \mathbb{N}_{\geq 2} \cup \{ \infty \}$,
with the convention that $\frac{\infty -1}{\infty} = 1$.
Such pairs are to be interpreted as geometric objects
interpolating between projective varieties and logarithmic pairs in the sense of Iitaka.
Differential forms on such pairs also interpolate between
regular Kähler differentials on projective varieties
and logarithmic differentials on logarithmic pairs.
More precisely, they are differential forms with ``fractional pole orders'',
e.g., on $(\mathbb{A}^{1},\frac{1}{2}\{ x = 0\})$ we would
want to consider something like $\frac{1}{\sqrt{x}}dx$.
A way to make this precise is by passing to a suitably ramified cover,
and this leads to the notion of \emph{adapted differentials}.
See \Cref{defn:adapteddifferentials} for the precise definition,
and \cite[\S 5]{cp19} for further motivation regarding adapted differentials.

Since we work with normal varieties,
we consider \emph{reflexive} differentials.
An important question regarding reflexive differentials
is whether they can be extended to a resolution of singularities.
More precisely, this means the following.
Given a normal variety $X$, a reflexive differential form $\sigma$ on $X$
can be understood as a regular Kähler differential form
defined only over the smooth locus $X_{\mathrm{reg}}$ of $X$.
Given a resolution of singularities $\pi \colon \tilde{X} \to X$,
which we assume for simplicity to be an isomorphism over $X_{\mathrm{reg}}$,
we identify $\tilde{X} \setminus \operatorname{Exc}(\pi)$ with $X_{\mathrm{reg}}$
and regard $\sigma$ as a regular Kähler differential form
on $\tilde{X} \setminus \operatorname{Exc}(\pi)$.
The differential form $\sigma$ extends to a rational
differential form $\tilde{\sigma}$ on the smooth variety $\tilde{X}$,
and this rational differential form may or may not have
poles along the divisor $\operatorname{Exc}(\pi)$.
In case it does not have any poles, i.e., in case $\tilde{\sigma}$ is
a regular Kähler differential on $\tilde{X}$, we say that $\sigma$ extends
to a regular differential form on the resolution of singularities.
This kind of extension property of reflexive differentials is
closely related to MMP singularities, cf.~\cite{gkkp11}.
Our goal here is to study this extension property in the case of
adapted differentials over klt $\mathcal{C}$-pairs.

\subsection{Main result}

Let $(X,\Delta)$ be a $\mathcal{C}$-pair.
Given a quasi-finite morphism $\gamma \colon Y \to X$ of normal varieties
of the same dimension such that $\gamma^{*}\Delta$ has integer coefficients,
which we refer to as an \emph{adapted morphism}, we can define
the \emph{sheaf of adapted differentials} $\Omega_{(X,\Delta,\gamma)}^{[1]}$ on $Y$;
see \Cref{defn:adapteddifferentials}.
If $\lfloor \Delta \rfloor = 0$, then $\Omega_{(X,\Delta,\gamma)}^{[1]} \subseteq \Omega_{Y}^{[1]}$
is a subsheaf of the sheaf of reflexive differentials on $Y$,
and we show the following:

\begin{thm}
  \label{thm:main}
  Let $(X,\Delta)$ be a $\mathcal{C}$-pair that is klt in the usual sense,
  and let $\gamma \colon Y \to X$ be an adapted morphism.
  Let $\pi \colon \tilde{Y} \to Y$ be a log resolution of singularities of $Y$.
  Then we can extend reflexive differential forms in $\Omega_{(X,\Delta,\gamma)}^{[1]}$
  to regular differential forms in $\Omega_{\tilde{Y}}^{1}$.
  That is, the pull-back of rational differential forms
  induces a morphism of $\mathscr{O}_{Y}$-modules
  \[ \Omega_{(X,\Delta,\gamma)}^{[1]} \to \pi_{*}\Omega_{\tilde{Y}}^{1}. \]
\end{thm}

\begin{rem}
  Such extension results are well-known for varieties with klt singularities,
  and are known to fail for varieties with worse than klt singularities;
  cf.~\cite{gkkp11}.
  A crucial point in \Cref{thm:main} is that, despite $(X,\Delta)$ being klt,
  the cover $Y$ may have worse than klt singularities.
  Therefore, we cannot hope to be able to extend \emph{any}
  reflexive differential form to the resolution of singularities.
  But the theorem ensures that we can extend \emph{adapted} reflexive differentials.
\end{rem}

\subsection{Relation to the literature}

Most of the contents of this note are extracted
from the author's PhD thesis \cite{thesis}.
In some cases we refer to \cite{thesis} for more detailed computations.
The results here generalize results of Kebekus--Rousseau in \cite{kr24}
to the singular setting.

\subsection{Acknowledgements}

I am very thankful to Stefan Kebekus for many useful discussions and ideas.
I would also like to thank Hsueh-Yung Lin for useful suggestions.

\section{Preliminaries}
\label{section:preliminaries}

\subsection{Notation and conventions}

We will follow notation and conventions from \cite{har77} and \cite{km98}.
For the most part, we will be interested in the algebraic setting.
However, some definitions and results need to be stated in the analytic setting,
to allow for analytic-local arguments.
Unless explicitly stated otherwise, open subsets will refer to Zariski-open subsets.

We follow \cite[II. \S 3]{del70} when it comes to logarithmic differentials,
and \cite[\S 2.E]{gkkp11} when it comes to reflexive sheaves.
When citing \href{https://stacks.math.columbia.edu/}{The Stacks Project},
we will use the format \cite[\href{https://stacks.math.columbia.edu/tag/0000}{0000}]{stacks-project}.
The four characters on the right represent the tag,
so the corresponding URL would be \href{https://stacks.math.columbia.edu/tag/0000}{https://stacks.math.columbia.edu/tag/0000}.

Given a Weil $\mathbb{Q}$-divisor $\Delta$ on a normal variety $X$,
and given a prime Weil divisor $D$ on $X$, we define $\operatorname{mult}_{D}(\Delta)$
to be the coefficient of $D$ in $\Delta$, i.e.,
the unique rational number $q \in \mathbb{Q}$ such that $D \not\subseteq \operatorname{Supp}(\Delta - qD)$.

\subsection{Pulling back reflexive differentials as rational differentials}
\label{subsection:pullback}

Let $X$ be a normal variety and let $p \in \mathbb{N}$.
A \emph{rational differential $p$-form} is a germ in $(\Omega_{X}^{p})_{\eta}$,
where $\eta$ denotes the generic point of $X$.
Equivalently, in the language of \cite[\href{https://stacks.math.columbia.edu/tag/01X1}{01X1}]{stacks-project},
it is a meromorphic section of $\Omega_{X}^{p}$, i.e., a global section
of $\mathscr{K}_{X}(\Omega_{X}^{p}) := \Omega_{X}^{p} \otimes_{\mathscr{O}_{X}} \mathscr{K}_{X}$,
where $\mathscr{K}_{X}$ is the constant sheaf with value $\mathscr{O}_{X,\eta}$.
Since Kähler differentials are torsion-free on smooth varieties,
specifying a rational differential $p$-form is the same as specifying
a rational section of $\Omega_{X_{\mathrm{reg}}}^{p}$, i.e., specifying a section
of $\Omega_{X_{\mathrm{reg}}}^{p}$ defined over some non-empty open subset of $X_{\mathrm{reg}}$.

Let now $f \colon X \to Y$ be a morphism of normal varieties.
If $f$ is dominant, then we have a well-defined pull-back
\[ f^{*} \colon (\Omega_{Y}^{p})_{\eta_{Y}} \to (\Omega_{X}^{p})_{\eta_{X}} \]
of rational differential $p$-forms.
Given a non-empty open subset $V \subseteq Y$,
the restriction map
\[ \Gamma(V,\Omega_{Y}^{[p]}) \to (\Omega_{Y}^{p})_{\eta_{Y}} \]
is injective, because reflexive differentials are torsion-free.
Therefore, a reflexive differential $p$-form $\sigma \in \Gamma(V,\Omega_{Y}^{[p]})$
can be regarded as a rational differential $p$-form,
and thus can be pulled back to a rational differential $p$-form
$f^{*}(\sigma)$ on $X$.
Denoting $U = f^{-1}(V)$, the restriction defines again
an injective map $\Gamma(U,\Omega_{X}^{[p]}) \to (\Omega_{X}^{p})_{\eta_{X}}$.
If $f^{*}(\sigma)$ is in the image of this map,
then it defines a uniquely determined reflexive differential $p$-form
over $U \subseteq X$, and we say that
\emph{the reflexive differential $\sigma$ on $V$ pulls back to a reflexive differential on $U$}.
If all reflexive differentials pull back to reflexive differentials,
then we obtain an $\mathscr{O}_{Y}$-module morphism
\[ f^{*} \colon \Omega_{Y}^{[p]} \to f_{*}\Omega_{X}^{[p]}, \]
and we say that this morphism is \emph{induced by the pull-back of rational differential forms}.
And similarly for other subsheaves of $\mathscr{K}_{Y}(\Omega_{Y}^{p})$ and $\mathscr{K}_{X}(\Omega_{X}^{p})$.

\subsection{$\mathcal{C}$-pairs}

\begin{defn}[{$\mathcal{C}$-pair, cf.~\cite[Definition 2.1]{cam11a}}]
  \label{defn:cpair}
  A \emph{$\mathcal{C}$-pair} $(X,\Delta)$ consists of
  a normal variety $X$ and a Weil $\mathbb{Q}$-divisor $\Delta$
  of the form
  \[ \Delta = \sum_{i \in I} \frac{m_{i}-1}{m_{i}}D_{i}, \]
  with $D_{i}$ pairwise distinct prime Weil divisors
  and with $m_{i} \in \mathbb{N}_{\geq 2} \cup \{ \infty \}$ for all $i \in I$,
  with the convention that $\frac{\infty -1}{\infty} = 1$.
\end{defn}

\begin{nota}
  \label{nota:boundary}
  Let $(X,\Delta)$ be a $\mathcal{C}$-pair.
  Then we denote $\Delta^{\mathrm{orb}} := \lceil \Delta \rceil - \Delta$.
  If $f \colon Y \to X$ is a dominant morphism such that
  we can pull back $\lfloor \Delta \rfloor$ along $f$ in a meaningful way,
  then we abuse notation and denote by $\lfloor \Delta \rfloor$
  the reduced Weil divisor $(f^{*}\lfloor \Delta \rfloor)_{\mathrm{red}}$ on $Y$.
\end{nota}

\begin{defnlm}[{Quotient $\mathcal{C}$-pair, cf.~\cite[\S 12]{kr24}}]
  \label{defnlm:quotient}
  Let $(X',\Delta')$ be a $\mathcal{C}$-pair and let $G \subseteq \operatorname{Aut}(X')$
  be a finite group acting on $X'$ such that $g(\Delta) = \Delta$ for all $g \in G$.
  Assume\footnote{
    This is automatically the case if $X'$ is quasi-projective.
  }
  that we can form the quotient of $X'$ by $G$
  and let $q \colon X' \to X$ be the quotient morphism.
  We define a $\mathcal{C}$-pair $(X,\Delta)$ as follows.
  For each prime Weil divisor $D$ on $X$,
  choose a prime Weil divisor $D'$ on $X'$
  contained in the support of $q^{*}D$.
  Then set
  \begin{equation}
    \label{eqn:mult}
    \operatorname{mult}_{D}(\Delta) :=
    \frac{\operatorname{mult}_{D'}(q^{*}D) \cdot \frac{1}{1 - \operatorname{mult}_{D'}(\Delta')}-1}
    {\operatorname{mult}_{D'}(q^{*}D) \cdot \frac{1}{1 - \operatorname{mult}_{D'}(\Delta')}},
  \end{equation}
  with the usual arithmetic conventions regarding $\infty$.
  The $\mathcal{C}$-pair $(X,\Delta)$ is called
  the \emph{quotient of the $\mathcal{C}$-pair $(X',\Delta')$ by $G$}.

  \begin{proof}
    We show that the resulting $\mathcal{C}$-pair is well-defined,
    i.e., that $\operatorname{mult}_{D}(\Delta)$ does not depend on the chosen $D'$.
    So let $D''$ be another prime Weil divisor on $X'$ contained in the support of $q^{*}D$.
    Then we can find some $g \in G$ such that $D' = g(D'')$.
    Moreover, since $g^{*}q^{*}D = q^{*}D$, $D'$ and $D''$ have the same coefficient in $q^{*}D$,
    i.e.,
    \[ \operatorname{mult}_{D''}(q^{*}D) = \operatorname{mult}_{D'}(q^{*}D). \]
    By assumption, we also have
    \[ \operatorname{mult}_{D''}(\Delta') = \operatorname{mult}_{D'}(\Delta'), \]
    hence the claim.
  \end{proof}
\end{defnlm}

\begin{rem}
  Quotient $\mathcal{C}$-pairs are a particular case of
  Campana's \emph{orbifold base} of a proper surjective morphism
  from a $\mathcal{C}$-pair to a normal variety,
  cf.~\cite[Définition 4.2]{cam11a}.
\end{rem} 

\begin{defn}[Singularities of $\mathcal{C}$-pairs]
  Let $(X,\Delta)$ be a $\mathcal{C}$-pair.
  We say that
    \begin{enumerate}
      \item the pair $(X,\Delta)$ is \emph{smooth}
      if $X$ is smooth and $\Delta$ is snc in the sense of \cite[Notation 0.4]{km98};
      \item the pair $(X,\Delta)$ has \emph{quotient singularities}
      if it is analytic-locally the quotient $\mathcal{C}$-pair
      of a smooth $\mathcal{C}$-pair by the action of a finite group;
      \item the pair $(X,\Delta)$ is \emph{klt}
      if $(X,\Delta)$ is klt in the sense of \cite[Definition 2.34]{km98}.
    \end{enumerate}
\end{defn}

\subsection{Adapted covers}

\begin{defn}[Quasi-cover]
  \label{defn:quasicover}
  A \emph{quasi-cover} is a quasi-finite morphism
  between normal varieties of the same dimension.
  A \emph{cover} is a finite quasi-cover, and a \emph{Galois cover}
  is a cover which is isomorphic to a geometric quotient
  by the action of a finite group.
\end{defn}

\begin{rem}
  \label{rem:quasicover}
  Quasi-covers are open and dominant, cf.~\cite[Lemma B.4]{thesis}.
  Covers are therefore surjective,
  but surjective quasi-covers need not be covers.
\end{rem}

\begin{defn}[Adapted morphism]
  \label{defn:adaptedcover}
  Let $(X,\Delta)$ be a $\mathcal{C}$-pair and
  let $\gamma \colon Y \to X$ be a quasi-cover.
  We say that
  \begin{enumerate}
    \item the morphism $\gamma$ is an \emph{adapted morphism}
    if $\gamma^{*}\Delta^{\mathrm{orb}}$ has integer coefficients;
    \item the morphism $\gamma$ is \emph{perfectly adapted}
    if $\gamma^{*}\Delta^{\mathrm{orb}}$ is a reduced divisor
    and $\gamma$ is quasi-\'etale over $X \setminus \operatorname{Supp}(\Delta)$.
  \end{enumerate}
\end{defn}

\begin{lm}
  \label{lm:quotientsingularitiesperfectlyadaptedcover}
  Let $(X,\Delta)$ be a $\mathcal{C}$-pair with quotient singularities.
  Then, for every point $x \in X$,
  there exists an analytic-open neighborhood $U$ of $x$ in $X$
  and a perfectly adapted cover $\gamma \colon V \to U$
  such that $V$ is smooth and $\gamma^{*}\lfloor \Delta \rfloor$ is snc.
  \begin{proof}
    We may assume that there is a smooth $\mathcal{C}$-pair $(X',\Delta')$
    such that $(X,\Delta)$ is the quotient
    of the $\mathcal{C}$-pair $(X',\Delta')$
    by some finite group $G$.
    Let $q \colon X' \to X$ denote the quotient morphism
    and let $z_{1}, \ldots, z_{n}$ be analytic-local coordinates
    around a point $z \in q^{-1}(x)$ such that $\Delta'$
    is locally given by $\sum_{i = 1}^{n} \frac{m_{i}-1}{m_{i}}\{z_{i} = 0 \}$ around $z$,
    for some $m_{1}, \ldots, m_{n} \in \mathbb{N}_{\geq 1} \cup \{ \infty \}$.
    Since we can shrink $X'$ and $X$ accordingly,
    we may assume that this is true globally.
    Let $I \subseteq \{ 1, \ldots, n \}$ be
    such that $m_{i} = \infty$ for all $i \in I$
    and $m_{i} \in \mathbb{N}_{\geq 1}$ for all $i \in \{1, \ldots, n \} \setminus I$,
    and define $a_{i} := 1$ if $i \in I$
    and $a_{i} := m_{i}$ if $i \in \{ 1, \ldots, n \} \setminus I$.
    Then we can write down a perfectly adapted cover $\delta \colon Y \to X'$ in these coordinates as
    \[ (y_{1}, \ldots, y_{n}) \mapsto (y_{1}^{a_{1}}, \ldots, y_{n}^{a_{n}}). \]
    A local computation shows that the composition $q \circ \delta \colon Y \to X$ is
    a perfectly adapted cover of $(X,\Delta)$,
    cf.~\cite[Corollary 1.55]{thesis}.
    And by construction $Y$ is smooth
    and $\gamma^{*}\lfloor \Delta \rfloor$ is snc.
  \end{proof}
\end{lm}

\begin{lm}
  \label{lm:qcartierperfectlyadaptedcover}
  Let $(X,\Delta)$ be a $\mathcal{C}$-pair
  such that every irreducible component of $\Delta^{\mathrm{orb}}$
  is $\mathbb{Q}$-Cartier.
  Then, for every point $x \in X$, there exists a Zariski-open
  neighborhood $U$ of $x$ in $X$ and a perfectly adapted cover $\gamma \colon V \to U$.
  \begin{proof}
    Let $N \in \mathbb{N}_{>0}$ be such that $ND_{i}$ is Cartier
    for every prime Weil divisor $D_{i}$ in the support of $\Delta^{\mathrm{orb}}$.
    Shrinking $X$ around $x$, we may assume that
    every such $ND_{i}$ is a principal effective Cartier divisor.
    If $\Delta^{\mathrm{orb}} = 0$,
    then we can take $\gamma = \operatorname{id}_{X}$.
    Therefore, we may assume that $\Delta^{\mathrm{orb}} = \sum_{i = 1}^{l} \frac{1}{m_{i}}D_{i}$
    with $l \geq 1$, and define $m := \operatorname{lcm}(m_{1}, \ldots, m_{l}) > 1$.
    Then $D := mN\Delta^{\mathrm{orb}}$ is a principal effective Cartier divisor,
    so we have an isomorphism $\mathscr{O}_{X}(D) \cong \mathscr{O}_{X}^{\otimes m}$,
    and we can apply the cyclic cover construction from \cite[Definition 2.52]{km98}.
    The resulting cover may be reducible,
    but after choosing an irreducible component,
    we obtain the desired perfectly adapted cover, cf.~\cite[Lemma 1.56]{thesis}.
  \end{proof}
\end{lm}

The covers constructed in \Cref{lm:quotientsingularitiesperfectlyadaptedcover,lm:qcartierperfectlyadaptedcover} need not be isomorphic.
The cover constructed in \Cref{lm:qcartierperfectlyadaptedcover} is always cyclic,
i.e., it is a Galois cover with cyclic automorphism group.
The cover constructed in \Cref{lm:quotientsingularitiesperfectlyadaptedcover} is also Galois,
but in this case the automorphism group of the cover need not be cyclic.
This is illustrated in the following example:

\begin{exmp}
  \label{exmp:covers}
  Let $X = \mathbb{A}^{2}$ with coordinates $x_{1}$ and $x_{2}$,
  and let $\Delta = \frac{1}{2}\{ x_{1} = 0 \} + \frac{1}{2}\{x_{2} = 0 \}$.
  Then, the cover produced by \Cref{lm:quotientsingularitiesperfectlyadaptedcover} is given by
  \[ \mathbb{A}^{2} \to \mathbb{A}^{2}, (z_{1},z_{2}) \mapsto (z_{1}^{2},z_{2}^{2}), \]
  so it is a Galois cover with non-cyclic Galois group $(\mathbb{Z}/2\mathbb{Z})^{2}$.
  On the other hand, the cover produced by \Cref{lm:qcartierperfectlyadaptedcover} is the cyclic cover given by
  \[ \{ (z_{1},z_{2},t) \in \mathbb{A}^{3} \mid t^{2} = z_{1}z_{2} \}
  \to \mathbb{A}^{2}, (z_{1},z_{2},t) \mapsto (z_{1},z_{2}). \]
  We refer to \cite[\S 4]{thesis} for more detailed computations
  and comparisons.
\end{exmp}

\begin{lm}
  \label{lm:kltcover}
  Let $(X,\Delta)$ be a $\mathcal{C}$-pair.
  If $\gamma \colon Y \to X$ is a perfectly adapted cover,
  then $(Y,0)$ is klt.
  \begin{proof}
    We have $\lfloor \Delta \rfloor = 0$, and we also have
    \[ K_{Y} = \gamma^{*}(K_{X}) + \operatorname{Ram}(\gamma) = \gamma^{*}(K_{X} + \Delta), \]
    so the claim follows from \cite[Proposition 5.20]{km98}.
  \end{proof}
\end{lm}

\begin{lm}
  \label{lm:kltcodimension2perfectlyadaptedcover}
  Let $(X,\Delta)$ be a klt $\mathcal{C}$-pair.
  Then there exists a closed subset $Z \subseteq X$ of codimension
  at least $3$ in $X$ with the following property.
  For every point $x \in X \setminus Z$ there exists
  an analytic-open neighborhood $U$ of $x$ in $X$
  and a perfectly adapted cover $\gamma \colon V \to U$
  such that $V$ is smooth.
  \begin{proof}
    After removing a closed subset of codimension $3$,
    we may assume that $X$ is $\mathbb{Q}$-factorial \cite[Proposition 9.1]{gkkp11}.
    By \Cref{lm:qcartierperfectlyadaptedcover},
    after shrinking $X$ further around $x$,
    we may assume that there exists a perfectly adapted cover $\gamma \colon Y \to X$.
    Then, \Cref{lm:kltcover} implies that $(Y,0)$ is klt.
    Since $Y$ has quotient singularities in codimension $2$ \cite[Proposition 9.3]{gkkp11}
    and $\gamma$ is finite and surjective,
    after removing another closed subset of codimension $3$ from $X$
    we may assume that every $y \in Y$ has an analytic-open neighborhood $V_{y}$
    which is a quasi-\'etale quotient of an affine space by a finite group action,
    say $q \colon \tilde{V} \to V_{y} = \tilde{V}/G$.
    By \cite[Theorem in p.~48]{gr84},
    we may assume that $\gamma|_{V_{y}} \colon V_{y} \to \gamma(V_{y})$
    is also a perfectly adapted cover.
    Therefore, since $q$ is quasi-\'etale and finite,
    the composition $\gamma \circ q \colon \tilde{V} \to \gamma(V_{y})$
    is a perfectly adapted cover with $\tilde{V}$ smooth.
  \end{proof}
\end{lm}

\subsection{Adapted differentials}

\begin{defn}[Adapted differentials]
  \label{defn:adapteddifferentials}
  Let $I$ be a set of indices such that $\Delta = \sum_{i \in I} \frac{m_{i} - 1}{m_{i}}D_{i}$
  with $D_{i}$ distinct prime Weil divisors.
  We denote by $I_{0} \subseteq I$ the subset of indices
  such that $m_{i} \in \mathbb{N}$.
  Let $\gamma \colon Y \to X$ be an adapted morphism
  and assume first that $(X,\Delta)$ and $(Y,\gamma^{*}\Delta)$
  are both snc pairs.
  Then $\gamma$ is flat,
  so $\gamma^{*}\mathscr{O}_{D_{i}} = \mathscr{O}_{\gamma^{*}D_{i}}$ for all $i \in I$,
  where $\gamma^{*}D_{i}$ is regarded here as a subscheme of $Y$
  with the appropriate non-reduced structure.
  Moreover, since $\gamma^{*}D_{i} \geq \gamma^{*}(\frac{1}{m_{i}}D_{i})$ for all $i \in I_{0}$,
  we have a quotient morphism
  \[ q \colon \oplus_{i \in I} \gamma^{*}\mathscr{O}_{D_{i}}
  \to \oplus_{i \in I_{0}} \mathscr{O}_{\gamma^{*}(\frac{1}{m_{i}}D_{i})} \]
  which is an epimorphism of $\mathscr{O}_{Y}$-modules.
  We define the $\mathscr{O}_{Y}$-module of \textit{adapted differential $1$-forms} as
  \[ \Omega_{(X,\Delta,\gamma)}^{1} :=
  \ker\left(\gamma^{*}\Omega_{X}^{1}(\log\lceil \Delta \rceil)
  \xrightarrow{\gamma^{*}(\operatorname{res})} \oplus_{i \in I} \gamma^{*} \mathscr{O}_{D_{i}}
  \xrightarrow{q} \oplus_{i \in I_{0}} \mathscr{O}_{\gamma^{*}(\frac{1}{m_{i}}D_{i})}\right). \]
  For each $p \in \mathbb{N}$,
  the $\mathscr{O}_{Y}$-module of \textit{adapted differential $p$-forms}
  is defined as
  \[ \Omega_{(X,\Delta,\gamma)}^{p} :=
  \bigwedge^{p} \Omega_{(X,\Delta,\gamma)}^{1}, \]
  which is locally free
  because $\Omega_{(X,\Delta,\gamma)}^{1}$ was locally free.

  In general, we consider the largest open subset $U \subseteq X$
  such that $(X,\Delta)$ is snc over $U$.
  This is a big open subset, i.e.,
  the codimension of its complement is at least $2$.
  Since $\gamma$ is a quasi-cover,
  $\gamma^{-1}(U)$ is a big open subset as well.
  We can then consider the largest open subset $V \subseteq \gamma^{-1}(U)$
  such that $(Y, \gamma^{*}\Delta)$ is snc over $V$,
  which is a big open subset of $Y$.
  We define the $\mathscr{O}_{Y}$-module of adapted differentials as
  \[ \Omega_{(X,\Delta,\gamma)}^{[p]} :=
  i_{*}\left(\Omega_{(U,\Delta|_{U},\gamma|_{V})}^{p}\right), \]
  where $i \colon V \to Y$ denotes the open immersion.
\end{defn}

\begin{rem}
  \label{rem:adaptedsnclocal}
  In \Cref{defn:adapteddifferentials},
  over suitable big open subsets,
  we may assume that $\gamma$ is given by
  \[ \gamma(y_{1}, \ldots, y_{n}) = (y_{1}^{a_{1}}, \ldots, y_{n}^{a_{n}}) \]
  for $a_{1}, \ldots, a_{n} \in \mathbb{N}_{\geq 1}$ with respect to
  suitable analytic-local coordinates for which $D_{i}$
  is locally given by $\{ x_{i} = 0 \}$ for all $i \in I$
  such that $D_{i}$ intersect the current analytic-local coordinate domain.
  The assumption that $\gamma$ is adapted translates
  into the fact that $m_{i} \mid a_{i}$ for all $i \in I_{0}$.
  Assume that $I_{0} = \{ 1, \ldots, l \}$,
  $I \setminus I_{0} = \{ l+1, \ldots, r \}$
  and $r < n$ for convenience and concreteness.
  Then, a set of local generators of
  the $\mathscr{O}_{Y}$-module $\Omega_{(X,\Delta,\gamma)}^{1}$ is given by
  \[ y_{1}^{\frac{a_{1}}{m_{1}}-1}dy_{1}, \ldots, y_{l}^{\frac{a_{l}}{m_{l}}-1}dy_{l},
  \frac{1}{y_{l+1}}dy_{l+1}, \ldots, \frac{1}{y_{r}}dy_{r}, y_{r+1}^{a_{r+1}-1}dy_{r+1},
  \ldots, y_{n}^{a_{n}-1}dy_{n}. \]
  See also \cite[\S 5.2]{cp19} or \cite[Remark 1.63]{thesis}.
\end{rem}

\begin{rem}[{cf.~\cite[\S 4]{kr24}}]
  \label{rem:reflexiveinclusions}
  Let $(X,\Delta)$ be a $\mathcal{C}$-pair
  and let $\gamma \colon Y \to X$ be an adapted morphism.
  Let $p \in \mathbb{N}$.
  Then $\Omega_{(X,\Delta,\gamma)}^{[p]}$ is reflexive and
  \[ \Omega_{(X,\Delta,\gamma)}^{[p]} \subseteq
  \Omega_{Y}^{[p]}(\log{\lfloor \Delta \rfloor}). \]
\end{rem}

\begin{lm}
  \label{lm:perfectlyadapted}
  Let $(X,\Delta)$ be a $\mathcal{C}$-pair
  and let $\gamma \colon Y \to X$ be a perfectly adapted morphism.
  Let $p \in \mathbb{N}$.
  Then we have
  \[ \Omega_{(X,\Delta,\gamma)}^{[p]} =
  \Omega_{Y}^{[p]}(\log{\lfloor \Delta \rfloor}). \]

  \begin{proof}
    By reflexivity, it suffices to show
    that the two sheaves agree over a big open subset.
    Therefore, we may assume that we are in the smooth and snc situation
    of \Cref{defn:adapteddifferentials} with $p = 1$.
    The claim follows then by direct computation,
    choosing appropriate analytic-local coordinates
    and using the local generators described in \Cref{rem:adaptedsnclocal}.
  \end{proof}
\end{lm}

\begin{lm}[{cf.~\cite[\S 4.4]{kr24}}]
  \label{lm:pullback}
  Let $\gamma_{1} \colon Y_{1} \to X$ and $\gamma_{2} \colon Y_{2} \to X$
  be two adapted morphisms and let $f \colon Y_{2} \to Y_{1}$ be a morphism over $X$.
  Then, for every $p \in \mathbb{N}$,
  the pull-back of rational differentials induces an isomorphism
  \[ f^{[*]}\Omega_{(X,\Delta,\gamma_{1})}^{[p]}
  \cong \Omega_{(X,\Delta,\gamma_{2})}^{[p]}. \]

  \begin{proof}
    By reflexivity of the sheaves involved,
    we may restrict our attention to suitable big open subsets.
    And since $f$ is a quasi-cover,
    the preimage of a big open subset under $f$ is again a big open subset.
    Therefore, we may assume that $X$, $Y_{1}$ and $Y_{2}$ are all smooth,
    and that $\Delta$, $\gamma_{1}^{*}\Delta$ and $\gamma_{2}^{*}\Delta$ are all snc.
    Then $\Omega_{(X,\Delta,\gamma_{1})}^{p}$ is a locally free subsheaf
    of $\Omega_{Y_{1}}^{p}(\log{\lfloor \Delta \rfloor})$ and $f$ is dominant,
    so the pull-back of rational differential forms induces an injective morphism
    \[ f^{*}\Omega_{(X,\Delta,\gamma_{1})}^{p}
    \hookrightarrow \Omega_{Y_{2}}^{p}(\log{\lfloor \Delta \rfloor}), \]
    cf.~\cite[\S 11.c]{iit82}.
    We claim that the image is $\Omega_{(X,\Delta,\gamma_{2})}^{p}$.
    It is enough to show this on the stalks and for $p = 1$,
    the morphism from a noetherian local ring to its completion
    is faithfully flat \cite[\href{https://stacks.math.columbia.edu/tag/00MC}{00MC}]{stacks-project},
    and the completion of a finitely generated module over said ring
    can be expressed as a tensor product with the completion
    of the ring \cite[Proposition 10.13]{am69},
    so by \cite[Proposition 3]{ser56} we may assume that
    we are in the analytic setting with $p=1$.
    The claim follows then by direct computation,
    choosing suitable analytic-local coordinates
    and using the local generators described in \Cref{rem:adaptedsnclocal}.
  \end{proof}
\end{lm}

\section{Extension properties of adapted differentials}
\label{section:extension}

As discussed in the introduction,
we are interested in extension of adapted reflexive differentials
to resolutions of singularities.
We first note the following:

\begin{lm}
  \label{lm:extensiondescent}
  Let $(X,\Delta)$ be a $\mathcal{C}$-pair.
  Let $\gamma_{1} \colon Y_{1} \to X$ be an adapted morphism
  and let $\pi_{1} \colon \tilde{Y}_{1} \to Y_{1}$ be
  a log resolution of singularities of $(Y_{1},\lfloor \Delta \rfloor)$.
  Suppose that we can find a commutative diagram
  \begin{center}
    \begin{tikzcd}
      \tilde{Y}_{2} \arrow{r}{\pi_{2}} \arrow[swap]{d}{\tilde{f}} & Y_{2} \arrow{d}{f} \\
      \tilde{Y}_{1} \arrow{r}{\pi_{1}} & Y_{1}
    \end{tikzcd}
  \end{center}
  such that $f$ is a quasi-cover,
  $\tilde{f}$ is surjective and generically finite,
  $\tilde{Y}_{2}$ is smooth and $\lfloor \Delta \rfloor$ is snc on $\tilde{Y}_{2}$.
  Suppose that, for $p \in \mathbb{N}$,
  the pull-back of rational differential forms
  induces an $\mathscr{O}_{Y_{2}}$-module morphism
  \[ \Omega_{(X,\Delta,\gamma_{1} \circ f)}^{[p]}
  \to (\pi_{2})_{*} \Omega_{\tilde{Y}_{2}}^{p}(\log{\lfloor \Delta \rfloor}). \]
  Then, the pull-back of rational differential forms
  induces an $\mathscr{O}_{Y_{1}}$-module morphism
  \[ \Omega_{(X,\Delta,\gamma_{1})}^{[p]}
  \to (\pi_{1})_{*} \Omega_{\tilde{Y}_{1}}^{p}(\log{\lfloor \Delta \rfloor}). \]
  \begin{proof}
    Let $V_{1} \subseteq Y_{1}$ be a non-empty open subset
    and let $\sigma \in \Omega_{(X,\Delta,\gamma_{1})}^{[p]}(V_{1})$ be an adapted differential.
    Let $V_{2} := f^{-1}(V_{1})$.
    By \Cref{lm:pullback}, the image of $\sigma$ under the composition
    \[ \Omega_{(X,\Delta,\gamma_{1})}^{[p]}(V_{1})
    \to f_{*}f^{*}\Omega_{(X,\Delta,\gamma_{1})}^{[p]}(V_{1})
    \to f_{*}f^{[*]}\Omega_{(X,\Delta,\gamma_{1})}^{[p]}(V_{1}) \]
    corresponds to a unique section $f^{*}(\sigma) \in \Omega_{(X,\Delta,\gamma_{1} \circ f)}^{[p]}(V_{2})$,
    which by assumption can be pulled back to a logarithmic differential
    \[ \pi_{2}^{*}f^{*}(\sigma) \in \Omega_{\tilde{Y}_{2}}^{p}(\log{\lfloor \Delta \rfloor})(\pi_{2}^{-1}(V_{2})). \]
    On the other hand, the pull-back $\pi_{1}^{*}(\sigma)$
    is a rational differential form on $\tilde{Y}_{1}$ such that
    \[ \tilde{f}^{*}\pi_{1}^{*}(\sigma) = \pi_{2}^{*}f^{*}(\sigma)
    \in \Omega_{\tilde{Y}_{2}}^{p}(\log{\lfloor \Delta \rfloor)}(\pi_{2}^{-1}(V_{2})). \]
    Since $\tilde{Y}_{1}$ is smooth
    and $\lfloor \Delta \rfloor \subseteq \tilde{Y}_{1}$ is snc,
    in order to show that $\pi_{1}^{*}(\sigma)$ is a logarithmic differential form,
    it suffices to show this over a big open subset of $\tilde{Y}_{1}$.
    Therefore, we may assume that $\tilde{f} \colon \tilde{Y}_{2} \to \tilde{Y}_{1}$ is a finite surjective morphism
    such that $\operatorname{Branch}(\tilde{f})$ and $\operatorname{Ram}(\tilde{f})$ are snc.

    Let $\tilde{Y}_{1}^{\circ} = \tilde{Y}_{1} \setminus \operatorname{Supp}(\lfloor \Delta \rfloor)$
    and $\tilde{Y}_{2}^{\circ} = \tilde{f}^{-1}(\tilde{Y}_{1}^{\circ})
    = \tilde{Y}_{2} \setminus \operatorname{Supp}(\lfloor \Delta \rfloor)$.
    Then we can apply \cite[Corollary 2.12.ii]{gkk10} to deduce
    that $\pi_{1}^{*}(\sigma)|_{\tilde{Y}_{1}^{\circ}}$ is a regular Kähler differential,
    because so is its pull-back along $\tilde{f}$.
    So it remains to show that $\pi_{1}^{*}(\sigma)$
    has at most logarithmic poles along $\lfloor \Delta \rfloor$.
    But $\operatorname{Supp}(\lfloor \Delta \rfloor)
    = \operatorname{Supp}(\tilde{f}^{-1}\lfloor \Delta \rfloor)$
    and $\tilde{f}^{*}\pi_{1}^{*}(\sigma)$ has at most logarithmic poles along $\lfloor \Delta \rfloor$,
    so we can apply \cite[Corollary 2.12.i]{gkk10} to conclude.
  \end{proof}
\end{lm}

We can now show the following extension result,
which will be used in the proof of \Cref{thm:main} later on:

\begin{lm}[{cf.~\cite[Construction 5.7]{kr24}}]
  \label{lm:niceextension}
  Let $(X,\Delta)$ be a $\mathcal{C}$-pair
  such that there exists a perfectly adapted cover $\gamma_{0} \colon Y_{0} \to X$
  with $Y_{0}$ smooth and $\gamma_{0}^{*}\lfloor \Delta \rfloor$ snc.
  Let $\gamma_{1} \colon Y_{1} \to X$ be an adapted morphism
  and let $\pi_{1} \colon \tilde{Y}_{1} \to Y_{1}$ be
  a log resolution of $(Y_{1},\lfloor \Delta \rfloor)$.
  Then, for any $p \in \mathbb{N}$,
  the pull-back of rational differential forms along $\pi_{1}$
  induces an $\mathscr{O}_{Y_{1}}$-module morphism
  \[ \Omega_{(X,\Delta,\gamma_{1})}^{[p]} \to (\pi_{1})_{*}\Omega_{\tilde{Y}_{1}}^{p}(\log{\lfloor \Delta\rfloor}). \]
  \begin{proof}
    Let $Y_{2}$ be an irreducible component
    of the normalization of the fiber product $Y_{1} \times_{X} Y_{0}$
    that surjects onto $Y_{1}$.
    Such a component exists because
    being finite and surjective is stable under base change.
    We obtain a commutative square
    \begin{center}
      \begin{tikzcd}
        Y_{2} \arrow[swap]{d}{f} \arrow{r}{\delta} & Y_{0} \arrow{d}{\gamma_{0}} \\
        Y_{1} \arrow{r}{\gamma_{1}} & X
      \end{tikzcd}
    \end{center}
    in which $f$ is a cover and $\delta$ is a quasi-cover.
    Let $\pi_{2} \colon \tilde{Y}_{2} \to Y_{2}$ be
    a log resolution of $(Y_{2},\lfloor \Delta \rfloor)$
    together with a surjective generically finite morphism $\tilde{f} \colon \tilde{Y}_{2} \to \tilde{Y}_{1}$
    fitting into a commutative diagram
    \begin{center}
      \begin{tikzcd}
        \tilde{Y}_{2} \arrow[swap]{d}{\tilde{f}} \arrow{r}{\pi_{2}} & Y_{2} \arrow{d}{f} \\
        \tilde{Y}_{1} \arrow{r}{\pi_{1}} & Y_{1},
      \end{tikzcd}
    \end{center}
    which can be constructed by first taking any log resolution
    and then resolving the indeterminacy of the corresponding rational map to $\tilde{Y}_{1}$.
    By \Cref{lm:extensiondescent},
    it suffices to show that the pull-back of rational differential forms
    induces an $\mathscr{O}_{Y_{2}}$-module morphism
    \[ \Omega_{(X,\Delta,\gamma_{0} \circ \delta)}^{[p]}
    \to (\pi_{2})_{*}\Omega_{\tilde{Y}_{2}}^{p}(\log{\lfloor \Delta \rfloor}). \]

    Since $\gamma_{0}$ is perfectly adapted,
    by \Cref{lm:perfectlyadapted}, we have
    \[ \Omega_{(X,\Delta,\gamma_{0})}^{[p]} = \Omega_{Y_{0}}^{[p]}(\log{\lfloor \Delta \rfloor}). \]
    And since $Y$ is smooth and $\gamma_{0}^{*}\lfloor \Delta \rfloor$ is snc,
    the latter is just the locally free sheaf $\Omega_{Y_{0}}^{p}(\log{\lfloor \Delta \rfloor})$
    of logarithmic differentials on $Y_{0}$.
    Its pull-back under $\delta$ is again locally free,
    hence reflexive, so \Cref{lm:pullback} yields an isomorphism
    \[ \delta^{*}\Omega_{Y_{0}}^{p}(\log{\lfloor \Delta \rfloor})
    \cong \Omega_{(X,\Delta,\gamma_{0} \circ \delta)}^{[p]}. \]
    Now, given an open subset $V_{2} \subseteq Y_{2}$
    and a section $\xi \in \delta^{*}\Omega_{Y_{0}}^{p}(\log{\lfloor \Delta \rfloor})(V_{2})$,
    we want to show that $\pi_{2}^{*}(\xi) \in \Omega_{\tilde{Y}_{2}}^{p}(\log{\lfloor \Delta \rfloor})(\pi_{2}^{-1}(V_{2}))$.
    This last statement is local on $\tilde{Y}_{2}$,
    hence on $Y_{2}$ and on $Y_{0}$,
    so we may assume that $Y_{2} = V_{2}$
    and that $Y_{0}$, $Y_{2}$ and $\tilde{Y}_{2}$ are all affine.
    Then we can write $\xi = \sum_{i \in I} \alpha_{i}\sigma_{i}$
    for some $\alpha_{i} \in \mathscr{O}_{Y_{2}}(Y_{2})$
    and some $\sigma_{i} \in \Omega_{Y_{0}}^{p}(\log{\lfloor \Delta \rfloor})(Y_{0})$,
    and we have
    \[ \pi_{2}^{*}(\xi) = \sum_{i \in I} (\alpha_{i} \circ \pi_{2})(\delta \circ \pi_{2})^{*}(\sigma_{i}). \]
    Since $\operatorname{Supp}(\lfloor \Delta \rfloor)
    = \operatorname{Supp}((\delta \circ \pi_{2})^{-1}\lfloor \Delta \rfloor)$
    and the $\sigma_{i}$ are logarithmic differentials on $(Y_{0},\lfloor \Delta \rfloor)$,
    we can apply \cite[Fact 2.9]{gkk10} to conclude that the pull-back $\pi_{2}^{*}(\xi)$
    is a logarithmic differential on $(\tilde{Y}_{2},\lfloor \Delta \rfloor)$.
  \end{proof}
\end{lm}

\begin{cor}
  \label{cor:quotientextension}
  Let $(X,\Delta)$ be a $\mathcal{C}$-pair
  with quotient singularities.
  Let $\gamma \colon Y \to X$ be an adapted morphism
  and let $\pi \colon \tilde{Y} \to Y$ be
  a log resolution of $(Y,\lfloor \Delta \rfloor)$.
  Then, for any $p \in \mathbb{N}$,
  the pull-back of rational differentials induces
  a morphism of $\mathscr{O}_{Y}$-modules
  \[ \Omega_{(X,\Delta,\gamma)}^{[p]} \to \pi_{*}\Omega_{\tilde{Y}}^{p}(\log{\lfloor \Delta \rfloor}). \]
  \begin{proof}
    The question is analytic-local on $Y$,
    so we may shrink $Y$ and $X$ and assume that there exists
    a perfectly adapted cover $\gamma_{0} \colon Y_{0} \to X$ such that $Y_{0}$ is smooth
    and $\gamma_{0}^{*}\lfloor \Delta \rfloor$ is snc,
    see \Cref{lm:quotientsingularitiesperfectlyadaptedcover}.
    The statement follows then from \Cref{lm:niceextension}.
  \end{proof}
\end{cor}

In order to apply \Cref{lm:niceextension}
to extend adapted $1$-forms
coming from klt $\mathcal{C}$-pairs,
we first need the following result,
inspired by Flenner's extension theorem \cite{fle88}
and by some of the ideas and techniques in \cite{ks21}:

\begin{lm}
  \label{lm:flenner}
  Let $X$ be a reduced complex space of pure dimension $n$.
  Let $\pi \colon \tilde{X} \to X$ be a resolution of singularities.
  Let $p \in \{ 0, \ldots, n \}$ and let $A \subseteq X$ be
  a closed complex subspace with $\operatorname{codim}_{X}(A) \geq p + 2$.
  Then, sections of $\pi_{*}\Omega_{\tilde{X}}^{p}$ extend uniquely across $A$.
  That is, if we denote by $j \colon X \setminus A \to X$ the inclusion,
  then the canonical map
  \[ \pi_{*}\Omega_{\tilde{X}}^{p} \to j_{*}j^{*}\pi_{*}\Omega_{\tilde{X}}^{p} \]
  is an isomorphism.
  \begin{proof}
    The statement is independent of the resolution,
    cf.~\cite[\S 1.3]{ks21},
    so we may assume that $\pi$ is a projective morphism
    which is an isomorphism over the smooth locus of $X$.
    Since the statement is local on $X$,
    we may also assume that $X$ is a closed complex subspace
    of an open ball $B \subseteq \mathbb{C}^{N}$ for some $N \in \mathbb{N}$.
    Denote by $i \colon X \to B$ the inclusion and let $\nu := i \circ \pi$ be the composition,
    which is then also a projective morphism.
    It suffices to show that sections of $\nu_{*}\Omega_{\tilde{X}}^{p}$
    extend uniquely across $A \subseteq B$.
    As explained in \cite[\S 2.4 and \S 8]{ks21},
    from Saito's version of the Decomposition Theorem for the projective morphism $\nu$,
    we obtain a (non-canonical) decomposition
    \[ \mathbf{R}\nu_{*}\Omega_{\tilde{X}}^{p} \cong K_{p} \oplus R_{p} \]
    into two cochain complexes $K_{p}, R_{p} \in \mathbf{D}_{\mathrm{coh}}^{\mathrm{b}}(\mathscr{O}_{B})$
    with the following properties:
    \begin{enumerate}
      \item\label{item:flenner1} The support of $R_{p}$ is contained in
      the singular locus $X_{\mathrm{sing}} \subseteq B$.
      In particular, since $\pi_{*}\Omega_{\tilde{X}}^{p}$ is torsion-free,
      we have $\mathscr{H}^{0}(R_{p}) = 0$.
      \item\label{item:flenner2} We have $\mathscr{H}^{k}(K_{p}) = 0$
      for all $k \geq n - p + 1$.
      \item\label{item:flenner3} We have
      $\mathbf{R}\mathscr{H}om_{\mathscr{O}_{B}}(K_{p},\omega_{B}^{\bullet}) \cong K_{n-p}[n]$.
    \end{enumerate}
    We want to apply \cite[Proposition 6.4]{ks21} to deduce
    that sections of $\nu_{*}\Omega_{\tilde{X}}^{p}
    = \mathscr{H}^{0}(\mathbf{R}\nu_{*}\Omega_{\tilde{X}}^{p})$
    extend uniquely across $A$.
    By property (\ref{item:flenner1}) it suffices to apply \cite[Proposition 6.4]{ks21}
    to the complex $K_{p}$,
    because $\mathscr{H}^{0}(\mathbf{R}\nu_{*}\Omega_{\tilde{X}}^{p}) = \mathscr{H}^{0}(K_{p})$.
    Since $\mathscr{H}^{j}(\mathbf{R}\nu_{*}\Omega_{\tilde{X}}^{p}) = 0$
    for all $j < 0$, the same is true for the complex $K_{p}$,
    so we may indeed apply \cite[Proposition 6.4]{ks21} to $K_{p}$.

    We want to show that
    \[ \dim(A \cap \operatorname{Supp}(R^{k}\mathscr{H}om_{\mathscr{O}_{B}}(K_{p},\omega_{B}^{\bullet})))
    \leq - (k + 2) \]
    for all $k \in \mathbb{Z}$,
    with the understanding that $\dim(\varnothing) \leq k$
    for any $k \in \mathbb{Z}$.
    
    Let first $k \in \mathbb{Z}$ be an integer
    such that $k \leq p - n$.
    Since $\dim(A) \leq n - p - 2$,
    we have
    \[ \dim(A \cap \operatorname{Supp}(R^{k}\mathscr{H}om_{\mathscr{O}_{B}}(K_{p},\omega_{B}^{\bullet})))
    \leq n - p - 2 \leq - (k + 2), \]
    so in this case the desired inequality holds.

    Let now $k \in \mathbb{Z}$ be an integer
    such that $k \geq p - n + 1$.
    It suffices to show that $R^{k}\mathscr{H}om_{\mathscr{O}_{B}}(K_{p},\omega_{B}^{\bullet}) = 0$,
    because then its support is empty
    and the dimension of the intersection with $A$ is
    smaller than any given integer.
    By property (\ref{item:flenner3}),
    we can rewrite this sheaf as $\mathscr{H}^{k}(K_{n-p}[n]) = \mathscr{H}^{k+n}(K_{n-p})$.
    And finally by property (\ref{item:flenner2}),
    we have $\mathscr{H}^{k+n}(K_{n-p}) = 0$,
    because $k \geq p - n + 1$.
  \end{proof}
\end{lm}

\begin{exmp}
  Let $X$ be a $3$-fold with an isolated singularity
  and let $\pi \colon \tilde{X} \to X$ be a resolution of singularities.
  Then, by \Cref{lm:flenner} and \cite[\href{https://stacks.math.columbia.edu/tag/0AY6}{0AY6}]{stacks-project},
  the sheaf $\pi_{*}\Omega_{\tilde{X}}^{1}$ is reflexive.
  Therefore, we can extend reflexive differential $1$-forms on $X$
  to the resolution of singularities.
  Indeed, a $1$-form defined on $X_{\mathrm{reg}}$ pulls back to
  a section of $\pi_{*}\Omega_{\tilde{X}}^{1}$ defined
  over the big open subset $X_{\mathrm{reg}}$,
  and by reflexivity we can extend this section
  to a global section of $\pi_{*}\Omega_{\tilde{X}}^{1}$,
  i.e., a globally defined regular Kähler differential $1$-form.
\end{exmp}

We now have all the necessary ingredients to prove \Cref{thm:main}:

\begin{proof}[Proof of \Cref{thm:main}]
  Let $(X,\Delta)$ be a klt $\mathcal{C}$-pair.
  Let $\gamma \colon Y \to X$ be an adapted morphism
  and let $\pi \colon \tilde{Y} \to Y$ be
  a log resolution of singularities of $Y$.
  We want to show that the pull-back of rational differentials
  induces a morphism of $\mathscr{O}_{Y}$-modules
  \[ \Omega_{(X,\Delta,\gamma)}^{[1]} \to \pi_{*}\Omega_{\tilde{Y}}^{1}. \]
  \Cref{lm:flenner} allows us to remove a closed subset
  of codimension at least $3$ from $Y$.
  Since $\gamma$ is a quasi-cover,
  this implies that we may also remove a closed subset
  of codimension at least $3$ from $X$.
  And since the statement is analytic-local on $X$,
  by \Cref{lm:kltcodimension2perfectlyadaptedcover},
  we may assume that there exists
  a perfectly adapted cover $\gamma_{0} \colon Y_{0} \to X$
  from a smooth variety $Y_{0}$.
  The result follows then from \Cref{lm:niceextension}.
\end{proof}

\section{Applications}
\label{section:applications}

An important application of the extension result \cite[Theorem 1.4]{gkkp11}
is the existence of a pull-back of reflexive differentials on klt spaces \cite[Theorem 4.3]{gkkp11},
from which many other applications follow, cf.~\cite[Part II]{gkkp11}.
In this section, we prove the analogous result
for adapted reflexive differential $1$-forms over klt $\mathcal{C}$-pairs.

In this case,
we cannot talk about pull-back induced by pull-back of rational differential forms
as in \Cref{subsection:pullback},
because the morphism along which we want to pull back
is not necessarily dominant.
Instead, we will simply talk about a morphism
which \emph{agrees with the pull-back of Kähler differentials wherever this makes sense},
and refer to the proof of the result
or the more detailed discussion in \cite{keb13}
for further clarification.

\begin{thm}
  Let $(X,\Delta)$ be a klt $\mathcal{C}$-pair.
  Let $\gamma \colon Y \to X$ be an adapted morphism
  and let $f \colon Z \to Y$ be a morphism from a normal variety
  whose image is not contained in the singular locus of $Y$.
  Then, there exists a morphism of $\mathscr{O}_{Y}$-modules
  \[ \Omega_{(X,\Delta,\gamma)}^{[1]} \to f_{*}\Omega_{Z}^{[1]} \]
  that agrees with the pull-back of Kähler differentials wherever this makes sense.
  \begin{proof}
    The proof is analogous to that of \cite[Theorem 4.3]{gkkp11}.
    Let $\pi_{Y} \colon \tilde{Y} \to Y$ be
    a log resolution of singularities
    which is an isomorphism over $Y_{\mathrm{reg}}$.
    Since $\pi_{Y}$ is surjective,
    there is some irreducible component $Z_{0}$
    of the fiber product $Z \times_{Y} \tilde{Y}$
    that surjects onto $Z$.
    Moreover, since $f(Z) \cap Y_{\mathrm{reg}} \neq \varnothing$
    and $\pi_{Y}$ is an isomorphism over $Y_{\mathrm{reg}}$,
    the natural projection $Z_{0} \to Z$ is birational.
    Let now $\tilde{Z} \to Z_{0}$ be
    a resolution of singularities,
    so that $\pi_{Z} \colon \tilde{Z} \to Z$ is also
    a resolution of singularities
    fitting into the following commutative diagram:
    \begin{center}
      \begin{tikzcd}
        \tilde{Z} \arrow[swap]{d}{\pi_{Z}} \arrow{r}{\tilde{f}} & \tilde{Y} \arrow{d}{\pi_{Y}} \\
        Z \arrow{r}{f} & Y.
      \end{tikzcd}
    \end{center}
    Let now $V \subseteq Y$ be an open subset
    and let $\sigma \in \Omega_{(X,\Delta,\gamma)}^{[1]}(V)$
    be an adapted differential $1$-form.
    By \Cref{rem:reflexiveinclusions} we can regard it as a section
    \[ \sigma \in \Omega_{Y}^{[1]}(V) = \Omega_{Y_{\mathrm{reg}}}^{1}(V^{\circ}), \]
    where $V^{\circ} := V \cap Y_{\mathrm{reg}}$.
    Let $W := f^{-1}(V)$,
    let $W^{\circ} := f^{-1}(V^{\circ}) \cap Z_{\mathrm{reg}}$
    and let $f_{V} \colon W^{\circ} \to V^{\circ}$ be the induced morphism.
    We want to define a section
    \[ f^{*}(\sigma) \in \Omega_{Z}^{[1]}(W) \]
    in a way that is compatible with the pull-back of Kähler differentials.
    If $W^{\circ} = \varnothing$,
    then set $f^{*}(\sigma) := 0$.
    Otherwise, the regular Kähler differential $\sigma \in \Omega_{Y_{\mathrm{reg}}}^{1}(V^{\circ})$ pulls back
    to a regular Kähler differential
    \[ f_{V}^{*}(\sigma) \in \Omega_{Z_{\mathrm{reg}}}^{1}(W^{\circ}), \]
    and since $W^{\circ} \neq \varnothing$ we can regard it
    as a rational differential form on $Z_{\mathrm{reg}}$.
    Therefore, we can pull it back as a rational differential to $\tilde{Z}$,
    and by commutativity of the diagram above we have an equality
    \[ \pi_{Z}^{*}(f_{V}^{*}(\sigma)) = \tilde{f}^{*}(\pi_{Y}^{*}(\sigma)) \]
    of rational differential forms.
    By \Cref{thm:main}, $\pi_{Y}^{*}(\sigma)$ defines
    a regular differential form on $\pi_{Y}^{-1}(V)$,
    so $\tilde{f}^{*}(\pi_{Y}^{*}(\sigma)) = \pi_{Z}^{*}(f_{V}^{*}(\sigma))$ defines a regular differential form on $\tilde{f}^{-1}\pi_{Y}^{-1}(V) = \pi_{Z}^{-1}(W)$.
    This implies that the rational differential form $f_{V}^{*}(\sigma)$
    on $Z_{\mathrm{reg}}$ defines a
    regular differential form on $W_{\mathrm{reg}} = W \cap Z_{\mathrm{reg}}$.
    So $f_{V}^{*}(\sigma)$ defines a reflexive differential form
    \[ f^{*}(\sigma) \in \Omega_{Z}^{[1]}(W). \]

    We check now that this defines a sheaf morphism.
    Let $V_{2} \subseteq V_{1} \subseteq Y$ be an inclusion of open subsets in $Y$
    and let $\sigma \in \Omega_{(X,\Delta,\gamma)}^{[1]}(V_{1})$ be an adapted differential $1$-form.
    Let $V_{i}^{\circ}$, $W_{i}$, $W_{i}^{\circ}$
    and $f_{V_{i}}$ be defined as above for $i \in \{1, 2\}$.
    The recipe above gives us reflexive differential forms
    \[ f^{*}(\sigma) \in \Omega_{Z}^{[1]}(W_{1})
    \text{ and } f^{*}(\sigma|_{V_{2}}) \in \Omega_{Z}^{[1]}(W_{2}), \]
    and we want to show that
    \[ f^{*}(\sigma)|_{W_{2}} = f^{*}(\sigma|_{V_{2}}). \]
    Suppose first that $W_{2}^{\circ} = \varnothing$,
    i.e., that $f^{-1}(V_{2}^{\circ}) = \varnothing$.
    Since $V_{2}^{\circ} = V_{2} \cap Y_{\mathrm{reg}}$
    and $W_{2} = f^{-1}(V_{2})$,
    this means that $f(W_{2}) \subseteq Y_{\mathrm{sing}}$.
    But $f$ is continuous and $Y_{\mathrm{sing}} \subseteq Y$ is closed,
    so $W_{2} \neq \varnothing$ would imply
    that $f(Z) \subseteq Y_{\mathrm{sing}}$,
    contradicting the assumption on $f$.
    Therefore, in this case, $W_{2} = \varnothing$
    and both sides of the equality are zero.

    Suppose now that $W_{2}^{\circ} \neq \varnothing$.
    In this case, both reflexive differentials
    are given by the same rational differential
    \[ f_{V_{1}}^{*}(\sigma)|_{W_{2}^{\circ}} = f_{V_{2}}^{*}(\sigma|_{V_{2}^{\circ}}), \]
    hence are equal by torsion-freeness of $\Omega_{Z}^{[1]}$.
  \end{proof}
\end{thm}

\printbibliography
\vfill

\end{document}